\documentclass[a4paper, 9pt, twosides]{amsart}

\usepackage{lmodern}
\usepackage{amssymb}
\usepackage{amsthm}
\usepackage{amsaddr}
\usepackage{mathtools}
\usepackage{MnSymbol}
\usepackage{graphicx}
\usepackage{caption}
\usepackage{subcaption}
\usepackage{url}

\usepackage{enumitem}
\setitemize{nolistsep}
\setenumerate{nolistsep}
\usepackage{todonotes}
\usepackage{algorithm,algorithmic}

\allowdisplaybreaks

\newtheorem{assump}{Assumption}
\newtheorem{example}{Example}
\newtheorem{remark}{Remark}
\newtheorem{theorem}{Theorem}

\newtheorem{defn}{Definition}
\newtheorem{prob}{Problem}

\newcommand{\norm}[1]{\left\lVert{#1}\right\rVert}
\newcommand{\abs}[1]{\left\lvert{#1}\right\rvert}

\newcommand{\nn}{\nonumber}

\newcommand{\pmat}[1]{\begin{pmatrix}#1\end{pmatrix}}

\newcommand{\R}{\mathbb{R}}
\newcommand{\N}{\mathbb{N}}
\renewcommand{\P}{\mathcal{P}}

\newcommand{\Sw}{\mathcal{S}}

\newcommand{\KL}{\mathcal{KL}}
\newcommand{\Kinfty}{\mathcal{K}_{\infty}}
\newcommand{\Lf}{\mathcal{L}}

\newcommand{\Nsw}{\mathrm{N}}

\newcommand{\T}{\mathrm{T}}

\newcommand{\stab}{\mathrm{S}}
\newcommand{\ustab}{\mathrm{U}}

\allowdisplaybreaks

\title[]{Yet another stability condition for switched nonlinear systems}
\author{Atreyee Kundu}
\address{Department of Electrical Engineering,\\Indian Institute of Technology Kharagpur, \\West Bengal - 721302, India,\\ E-mail: atreyee@ee.iitkgp.ac.in}

\keywords{Switched Systems, Input/Output-to-State Stability, Multiple-IOSS Lyapunov-like Functions}

\date{\today}

\begin{document}

	\begin{abstract}
         We report a new class of switching signals that preserves input/output-to-state stability (IOSS) of continuous-time switched nonlinear systems under pre-specified restrictions on admissible switches between the subsystems and admissible dwell times on the subsystems. The primary apparatus for our analysis is multiple Lyapunov-like functions. Input-to-state stability (ISS) and global asymptotic stability (GAS) of switched systems under pre-specified restrictions on switching signals fall as special cases of our results when no outputs (resp., also inputs) are considered.
    \end{abstract}

    \maketitle
\section{Introduction}
\label{s:intro}
	Switched systems \cite[Section 1.1.2]{Liberzon2003} find wide applications in power systems, power electronics, aircraft control, networked control, etc.. In this paper we study input/output-to-state stability (IOSS) of continuous-time switched nonlinear systems. IOSS implies that irrespective of the initial state, if the inputs and the observed outputs are small, then the state of the system will become small eventually. This property also plays a key role in the design of state-norm estimators for switched systems. 

    IOSS of switched differential inclusions under arbitrary switching was studied in \cite{Mancilla2005}. The analysis relies on the existence of a common IOSS Lyapunov function. In \cite{Sanfelice_2010} IOSS of a class of hybrid systems that admits a Lyapunov function satisfying an IOSS relation both along the flow and during the jumps, was addressed. IOSS of switched systems under average dwell time switching \cite[Chapter 3]{Liberzon2003} was studied in \cite{Muller2012}. It was shown that if a switching signal obeys an average dwell time property along with a constrained duration of activation for the unstable subsystems, then IOSS of the resulting switched system is preserved. IOSS of impulsive switched systems is addressed in \cite{Li_2018}. Stability conditions for both stabilizing and destabilizing impulses are presented by employing Lyapunov functions and average dwell time switching conditions.


    In this paper we focus on IOSS of switched systems whose switching signals obey pre-specified restrictions on admissible switches between the subsystems and admissible dwell times on the subsystems. Such restrictions on switching signals often arise in many engineering systems, see \cite[Remark 1]{abc} for examples. 
    
    Earlier in \cite{abc,def} we addressed the regime of restricted switching. In \cite{abc} we presented an algorithm to design \emph{a stabilizing switching signal} that obeys given restrictions. The said design involves finding a class of cycles on a weighted directed graph representation of a switched system that satisfies certain conditions involving multiple Lyapunov-like functions \cite{Branicky1998} and the admissible dwell times. In \cite{def} we focussed on identifying \emph{a class of stabilizing switching signals} that obeys given restrictions. We considered common admissible minimum and maximum dwell times on all the subsystems, and our stability conditions relied on two factors: (i) non-consecutive activation of unstable subsystems and (ii) an inequality involving a set of scalars computed from Lyapunov-like functions corresponding to the individual subsystems and an admissible choice of dwell times on the subsystems. In the current work we continue our study of \emph{classes of} stabilizing switching signals that obey pre-specified restrictions. We accommodate (possibly different) admissible dwell times on the subsystems and relax the restriction of non-consecutive activation of unstable subsystems.  
	
	Our specific contribution is the following: Given a family of systems, possibly containing unstable dynamics, a set of pre-specified restrictions on admissible switches between the subsystems and admissible dwell times on the subsystems, we present a new class of switching signals that obeys the given restrictions and preserves IOSS of the resulting switched system. 
We show that if a certain inequality involving (a) a set of scalars in a given range, (b) a set of scalars computed from Lyapunov-like functions corresponding to the individual subsystems and the admissible transitions between them, and (c) the set of admissible minimum and maximum dwell times on the subsystems, is satisfied, then every switching signal that obeys the given restrictions with the  frequency of activation of the subsystems and the frequency of switches between the pairs of subsystems governed by functions of the set of scalars mentioned in (a), preserves IOSS of a switched system. In the absence of outputs (resp., also inputs), our class of switching signals ensures input-to-state stability (ISS) (resp., global asymptotic stability (GAS)) of a switched system. 
To the best of our knowledge, this is the first instance in the literature where \emph{a class of} switching signals that obeys pre-specified restrictions on admissible switches between the subsystems and admissible (possibly different) dwell times on the subsystems and preserves IOSS of the resulting switched systems, is characterized without a restriction on consecutive activation of unstable subsystems.

	  {\bf Notation}. \(\R\) is the set of real numbers, \(\norm{\cdot}\) is the Euclidean norm, and for any interval \(I\subseteq[0,+\infty[\), \(\norm{\cdot}_{I}\) is the essential supremum norm of a map from \(I\) into some Euclidean space. For a real number \(a\), we denote by \(\Bigl\lfloor a\Bigr\rfloor\) the biggest integer that is smaller than or equal to \(a\). For a finite set \(A\), we denote its cardinality by \(\abs{A}\).
\section{Problem statement and preliminaries}
\label{s:prob_stat}
	 We consider a family of continuous-time systems
    \begin{align}
    \label{e:family}
    \begin{aligned}
        \dot{x}(t) &= f_{p}(x(t),v(t)),\\
        y(t) &= h_{p}(x(t)),
    \end{aligned}
    \:\:x(0) = x_{0}\:\text{(given)},\:p\in\P,\:t\geq 0,
    \end{align}
    where \(x(t)\in\R^{d}\), \(v(t)\in\R^{m}\) and \(y(t)\in\R^{p}\) are the vectors of states, inputs and outputs at time \(t\), respectively, and \(\P = \{1,2,\ldots,N\}\) is an index set. Let \(\sigma:[0,+\infty[\to\P\) be a \emph{switching signal} --- it is a piecewise constant function that selects at each time \(t\), the index of the {active subsystem}, i.e., the system from the family \eqref{e:family} that is currently being followed. By convention, \(\sigma\) is assumed to be continuous from right and having limits from the left everywhere. Let \(\Sw\) denote the set of all switching signals. The \emph{switched system} generated by a family of systems \eqref{e:family} and a switching signal \(\sigma\) is given by
    \begin{align}
    \label{e:swsys}
    \begin{aligned}
        \dot{x}(t) &= f_{\sigma(t)}(x(t),v(t)),\\
        y(t) &= h_{\sigma(t)}(x(t)),
    \end{aligned}
    \:\:x(0) = x_{0}\:\text{(given)},\:t\geq 0.
    \end{align}
    We assume that for each \(p\in\P\), \(f_{p}\) is locally Lipschitz, \(f_{p}(0,0) = 0\) and \(h_{p}\) is continuous, \(h_{p}(0) = 0\); the exogenous inputs are Lebesgue measurable and essentially bounded. Thus, a solution to the switched system \eqref{e:swsys} exists in Carath\'eodory sense for some non-trivial time interval containing \(0\) \cite[Chapter 2]{Filippov1988}.

       \begin{defn}{\cite[Appendix A.6]{Liberzon2003}}
    	\label{d:ioss-cont}
	\rm{
        		The switched system \eqref{e:swsys} is \emph{input/output-to-state stable (IOSS)} under a switching signal \(\sigma\in\mathcal{S}\) if there exist class \(\mathcal{K}_{\infty}\) functions \(\alpha\), \(\chi_{1}\), \(\chi_{2}\) and a class \(\mathcal{KL}\) function \(\beta\) such that for all inputs \(v\) and initial states \(x_{0}\), we have 
        		\begin{align}
        		\label{e:ioss-cont}
            		\alpha(\norm{x(t)}) \leq \beta(\norm{x_{0}},t) + \chi_{1}(\norm{v}_{[0,t]}) + \chi_{2}(\norm{y}_{[0,t]})
        		\end{align}
        for all \(t\geq 0\).  
        If \(\chi_{2}\equiv 0\), then \eqref{e:ioss-cont} reduces to input-to-state stability (ISS) of \eqref{e:swsys}, and if also \(v\equiv 0\), then \eqref{e:ioss-cont} reduces to global asymptotic stability (GAS) of \eqref{e:swsys}.
      }
    \end{defn}

    Let \(\P_{S}\) and \(\P_{U}\) denote the sets of indices of IOSS and non-IOSS subsystems, respectively, \(\P = \P_{S}\sqcup\P_{U}\). Let \(E(\P)\) denote the set all pairs \((p,q)\) such that it is allowed to switch from subsystem \(p\) to subsystem \(q\), \(p,q\in\P\), \(p\neq q\). We let \(0=:\tau_{0}<\tau_{1}<\cdots\) be the \emph{switching instants}; these are the points in time where \(\sigma\) jumps. Let \(\delta_p\) and \(\Delta_p\) denote the admissible minimum and maximum dwell time on subsystem \(p\in\P\), respectively.
    \begin{defn}
    \label{d:adm-sw}
    \rm{
        A switching signal \(\sigma\in\Sw\) is called \emph{admissible} if it obeys the following two conditions:
        \begin{align}
        \label{res:adm_sw} \bigl(\sigma(\tau_i),\sigma(\tau_{i+1})\bigr)\in E(\P),\:\:i=0,1,\ldots,\\
        \intertext{and}
	    \label{res:adm_dw}
            \tau_{i+1}-\tau_{i} \in [\delta_{\sigma(\tau_i)},\Delta_{\sigma(\tau_i)}],\:\:
        i=0,1,2,\ldots.
        \end{align}
    }
    \end{defn}
     Let \(\Sw_{\mathcal{R}}\) denote the set of all admissible switching signals \(\sigma\). We will solve the following problem:
     \begin{prob}
     \label{prob:main}
        Given a family of systems \eqref{e:family}, the set of admissible switches between the subsystems, \(E(\P)\), and the admissible minimum and maximum dwell times, \(\delta_p\) and \(\Delta_p\), on the subsystems \(p\in\P\), identify a class of admissible switching signals \(\tilde{\Sw}_{\mathcal{R}}\subseteq\Sw_{\mathcal{R}}\) that ensures IOSS of the switched system \eqref{e:swsys}.
     \end{prob}
%
		 In short, our aim is to identify, if exists, a subset of \(\Sw_\mathcal{R}\) whose elements are stabilizing.
	 We will employ the temporal behaviour of Lyapunov-like functions for the systems in family \eqref{e:family} along the corresponding system trajectories.
     \begin{assump}
    \label{assump:key1}
        There exist class \(\Kinfty\) functions \(\underline{\alpha}\), \(\overline{\alpha}\), \(\gamma_{1}\), \(\gamma_{2}\), continuously differentiable functions \(V_{p}:\R^{d}\to[0,+\infty[\), \(p\in\P\), and constants \(\R\ni\lambda_p > 0\), \(p\in\P_S\) and \(\R\ni\lambda_p < 0\), \(p\in\P_U\), such that for all \(\xi\in\R^{d}\) and \(\eta\in\R^{m}\), the following hold:
        \begin{align}
        \label{e:key_ineq1}
            \underline{\alpha}(\norm{\xi})&\leq V_{p}(\xi)\leq\overline{\alpha}(\norm{\xi}),\:p\in\P,\\
            \intertext{and}
        \label{e:key_ineq2}
            \frac{\partial V_p}{\partial \xi}f_p(\xi,\eta)&\leq-\lambda_{p}V_{p}(\xi) + \gamma_{1}\bigl(\norm{\eta}\bigr) + \gamma_{2}\biggl(\norm{h_{p}(\xi)}\biggr),\:p\in\P.
         \end{align}
    \end{assump}
     \begin{assump}
    \label{assump:key2}
        There exist \(\mu_{pq} {> 0}\), \((p,q)\in E(\P)\) such that the IOSS-Lyapunov-like functions are related as follows:
        \begin{align}
        \label{e:key_ineq3}
            V_{q}(\xi)\leq\mu_{pq} V_{p}(\xi)\:\:\text{for all}\:\:\xi\in\R^{d}.
        \end{align}
    \end{assump}

     The functions \(V_{p}\), \(p\in\P\), are called the \emph{(multiple) IOSS-Lyapunov-like functions}. Condition \eqref{e:key_ineq2} is equivalent to the IOSS property for IOSS subsystems \cite{Krichman2001,Sontag1995} and the unboundedness observability property for the non-IOSS subsystems \cite{Angeli1999,Sontag1997}. The scalars \(\lambda_{p}\), \(p\in\P_S\) (resp., \(\lambda_p\), \(p\in\P_U\)) provide a quantitative measure of stability (resp., instability) of the subsystems \(p\in\P\). Condition \eqref{e:key_ineq3} restricts the class of Lyapunov-like functions to be linearly comparable.
     Given a family of systems \eqref{e:family}, the choice of IOSS-Lyapunov-like functions, \(V_p\), corresponding to the subsystems, \(p\in \P\), and consequently, the scalars, \(\lambda_p\), \(p\in\P\) and \(\mu_{pq}\), \((p,q)\in E(\P)\) is not unique. Let us consider that the functions \(V_p\), \(p\in \P\) and their corresponding scalars \(\lambda_p\), \(p\in\P\) and \(\mu_{pq}\), \((p,q)\in E(\P)\) are ``given''. We define
     \(E_{-}(\P) = \{(p,q)\in E(\P)\:|\:\ln\mu_{pq} < 0\}\) and \(E_{+}(\P) = \{(p,q)\in E(\P)\:|\:\ln\mu_{pq} > 0\}\).
    
    Fix an interval \(]s,t]\subseteq[0,+\infty[\) of time.
        (a) Let \(\Nsw(s,t)\) denote the total number of switches on \(]s,t]\).
        (b) We let \(\Nsw_p(s,t)\) denote the total number of switches on \(]s,t]\) where the subsystem \(p\) is activated.
        (c) Let \(\Nsw_{pq}(s,t)\) denote the total number of switches on \(]s,t]\) where a switch from subsystem \(p\) to subsystem \(q\) has occurred.
        (d) Let \(\T_{p}(s,t)\) denote the total duration of activation of the subsystem \(p\) on \(]s,t]\).
    We suppress the dependence of the above quantities on \(\sigma\) for notational simplicity. Further,
      		 the quantity \(\frac{\Nsw_p(s,t)}{\Nsw(s,t)}\) gives the \emph{frequency of activation of subsystem \(p\in\P\)} on \(]s,t]\), and
		 the quantity \(\frac{\Nsw_{pq}(s,t)}{\Nsw(s,t)}\) gives the \emph{frequency of switches from subsystem \(p\) to subsystem \(q\) on \(]s,t]\), \((p,q)\in E(\P)\)}.


\section{Main result}
\label{s:mainres}	
	\begin{assump}
	\label{assump:key3}
		There exist \(\rho_{p}^{S}\in[0,1[\), \(p\in\P_S\), \(\rho_{p}^{U}\in[0,1[\), \(p\in\P_U\) and \(\rho^{-}_{pq}\in[0,1[\), \((p,q)\in E_{-}(\P)\), \(\rho^{+}_{pq}\in[0,1[\), \((p,q)\in E_{+}(\P)\) satisfying \(\displaystyle{\sum_{p\in\P_S}\rho^{S}_p\leq 1}\), \(\displaystyle{\sum_{p\in\P_U}\rho^{U}_p < 1}\), \(\displaystyle{\sum_{(p,q)\in E_{-}(\P)}\rho^{-}_{pq}\leq 1}\), \(\displaystyle{\sum_{(p,q)\in E_{+}(\P)}\rho^{+}_{pq}\leq 1}\) such that
		\begin{align}
     	\label{e:key_ineq4}
			&-\frac{1}{\Delta_{\max}}\Biggl(\sum_{p\in\P_S}\abs{\lambda_p}\rho^{S}_{p}\delta_p+\sum_{(p,q)\in E_{-}(\P)}\abs{\ln\mu_{pq}}\rho^{-}_{pq}\Biggr)\nonumber\\
			&+\frac{1}{\delta_{\min}}\Biggl(\sum_{p\in\P_U}\abs{\lambda_p}\rho^{U}_{p}\Delta_p+\sum_{(p,q)\in E_{+}(\P)}\abs{\ln\mu_{pq}}\rho^{+}_{pq}\Biggr) < 0,
		\end{align}
		where \(\displaystyle{\delta_{\min} = \min_{p\in\P}\delta_p}\), \(\displaystyle{\Delta_{\max} = \max_{p\in\P}\Delta_p}\).
	\end{assump}
	
	Assumption \ref{assump:key3} requires the existence of four sets of scalars: (i) \(\rho_{p}^{S}\), \(p\in\P_S\), (ii) \(\rho_{p}^{U}\), \(p\in\P_U\), (iii) \(\rho^{-}_{pq}\), \((p,q)\in E_{-}(\P)\) and (iv) \(\rho^{+}_{pq}\), \((p,q)\in E_{+}(\P)\) that together with (a) the scalars \(\lambda_p\), \(p\in\P\) and \(\mu_{pq}\), \((p,q)\in E(\P)\) obtained from the IOSS-Lyapunov-like functions of the individual subsystems \(p\in\P\) and (b) the admissible minimum and maximum dwell times, \(\delta_p\) and \(\Delta_p\), on the subsystems \(p\in\P\), respectively, satisfy the inequality \eqref{e:key_ineq4}. We will use the scalars \(\rho_{p}^{S}\), \(p\in\P_S\), \(\rho_{p}^{U}\), \(p\in\P_U\), \(\rho^{-}_{pq}\), \((p,q)\in E_{-}(\P)\) and \(\rho^{+}_{pq}\), \((p,q)\in E_{+}(\P)\) to bound the frequency of activation of the stable subsystems, the frequency of activation of the unstable subsystems and the frequency of switches between the pairs of subsystems, respectively.
	
	\begin{defn}
	\label{d:adm-sw}
		Let \({\mathcal{S}}_{\mathcal{R}}'\) be the set of all switching signals, \(\sigma\in\mathcal{S}_{\mathcal{R}}\), that obey the following properties on every interval \(]s,t]\subseteq[0,+\infty[\) of time:\footnote{Notice that without the upper bounds on the sums of \(\rho^{S}_p\),\(p\in\P_S\), \(\rho^{U}_p\),\(p\in\P_u\), \(\rho^{+}_p\),\((p,q)\in E_{+}(\P)\), 
\(\rho^{-}_p\),\((p,q)\in E_{-}(\P)\) conditions \eqref{e:key_ineq4b}-\eqref{e:key_ineq4c} cannot be ensured. Also, the strict upper bound on the sum of \(\rho^{U}_p\), \(p\in\P_U\) ensures the activation of stable subsystems.}
		\begin{align}
            \label{e:key_ineq4a}&\biggl\lfloor\frac{t-s}{\Delta_{\max}}\biggr\rfloor\leq\Nsw(s,t)\leq\biggl\lfloor\frac{t-s}{\delta_{\min}}\biggr\rfloor,\\
			\label{e:key_ineq5}&\Nsw_p(s,t)\geq\Bigl\lfloor\rho^{S}_{p}\Nsw(s,t)\Bigr\rfloor,\:\:p\in\P_S,\\
			\label{e:key_ineq6}&\Nsw_p(s,t)\leq\tilde{\rho}^{\mathrm{U}}_{p}+\Bigl\lfloor\rho^{U}_{p}\Nsw(s,t)\Bigr\rfloor,\:\:p\in\P_U,\\
			\label{e:key_ineq7}&\Nsw_{pq}(s,t)\geq\Bigl\lfloor\rho^{-}_{pq}\Nsw(s,t)\Bigr\rfloor,\:\:(p,q)\in E_{-}(\P),\\
			\label{e:key_ineq8}&\Nsw_{pq}(s,t)\leq\tilde{\rho}^{+}_{pq}+\Bigl\lfloor\rho^{+}_{pq}\Nsw(s,t)\Bigr\rfloor,\:\:(p,q)\in E_{+}(\P),\\
            \label{e:key_ineq4b}&\sum_{p\in\P_S}\Nsw_{p}(s,t)+\sum_{p\in\P_U}\Nsw_{p}(s,t)=\Nsw(s,t),\\
            \label{e:key_ineq4c}&\sum_{(p,q)\in E_{-}(\P)}\Nsw_{pq}(s,t)+\sum_{(p,q)\in E_{+}(\P)}\Nsw_{pq}=\Nsw(s,t).
		\end{align}
        Here \(\tilde{\rho}^{\mathrm{U}}_{p}\), \(p\in\P_U\) and \(\tilde{\rho}^{+}_{pq}\), \((p,q)\in E_{+}(\P)\) are positive integers, and the scalars \(\rho^{\stab}_{p}\), \(p\in\P_S\), \(\rho^{\ustab}_{p}\), \(p\in\P_U\), \(\rho^{+}_{pq}\), \((p,q)\in E_{+}(\P)\), \(\rho^{-}_{pq}\), \((p,q)\in E_{-}(\P)\) are as described in Assumption \ref{assump:key3}.
	\end{defn}

    The elements of \(\mathcal{S}_{\mathcal{R}}'\) obey, on every interval of time, restrictions on (a) the total number of switches (condition \eqref{e:key_ineq4a}), (b) the minimum number of switches when stable subsystems \(p\) are activated (condition \eqref{e:key_ineq5}), (c) the maximum number of switches when unstable subsystems \(p\) are activated (condition \eqref{e:key_ineq6}), (d) the minimum number of switches from subsystems \(p\) to subsystems \(q\) such that \((p,q)\in E_{-}(\P)\) (condition \eqref{e:key_ineq7}), (e) the maximum number of switches from subsystems \(p\) to subsystems \(q\) such that \((p,q)\in E_{+}(\P)\) (condition \eqref{e:key_ineq8}), (f) the total number of switches where a stable (resp., unstable) subsystem is activated (condition \eqref{e:key_ineq4b}), and (g) the total number of switches between the admissible pairs of subsystems (condition \eqref{e:key_ineq4c}). Since condition \eqref{res:adm_dw} implies condition \eqref{e:key_ineq4a} and no property of \(\sigma\in\mathcal{S}_{\mathcal{R}}'\) violates condition \eqref{res:adm_sw}, it follows that
        {\it
            \(\mathcal{S}_{\mathcal{R}}'\subseteq\mathcal{S}_{\mathcal{R}}\).
        }
%
       
        We note that: 
        		\eqref{e:key_ineq5} provides a lower bound on the frequency of activation of the subsystems in \(\P_S\), 
		\(\displaystyle{\frac{\Nsw_{p}(s,t)}{\Nsw(s,t)}\geq\rho^{\stab}_{p}-1}\), \(p\in\P_S\),
		\eqref{e:key_ineq6} provides an upper bound on the frequency of activation of the subsystems in \(\P_U\),
		\(\displaystyle{\frac{\Nsw_{p}(s,t)-\tilde{\rho}^{\mathrm{U}}_{p}}{\Nsw(s,t)}\leq\rho^{\ustab}_{p}}\), \(p\in\P_U\),
		\eqref{e:key_ineq7} provides a lower bound on the frequency of switches from subsystems \(p\) to subsystems \(q\) such that \((p,q)\in E_{-}(\P)\),
		\(\displaystyle{\frac{\Nsw_{pq}(s,t)}{\Nsw(s,t)}\geq\rho^{-}_{pq}-1}\), \((p,q)\in E_{-}(\P)\), and
		\eqref{e:key_ineq8} provides an upper bound on the frequency of switches from subsystems \(p\) to subsystems \(q\) such that \((p,q)\in E_{+}(\P)\), 
		\(\displaystyle{\frac{\Nsw_{pq}(s,t)-\tilde{\rho}^{+}_{pq}}{\Nsw(s,t)}\leq\rho^{+}_{pq}}\), \((p,q)\in E_{+}(\P)\).
    \begin{remark}
    \label{rem:sw-prop2}
    \rm{
        Without the integers \(\tilde{\rho}^{\mathrm{U}}_{p}\), \(p\in\P_U\) and \(\tilde{\rho}^{+}_{pq}\), \((p,q)\in E_{+}(\P)\), we can never activate subsystems \(p\in\P_U\) and/or switch from subsystems \(p\) to subsystems \(q\) with \((p,q)\in E_{+}(\P)\). Indeed, consider an interval \(]s,t]\), where \(\Nsw(s,t) = 10\) is an admissible number of switches. Suppose that \(\tilde{\rho}^{\mathrm{U}}_{p}=0\), \(p\in\P_U\). Let us try to activate subsystem \(p'\in\P_U\) at a time \(s'\in]s,t]\). However, on an interval \(]s'-\varepsilon,s'+\varepsilon]\), \(\varepsilon > 0\) (small enough) such that \(\Nsw(s'-\varepsilon,s'+\varepsilon) = 0\), we are not allowed to activate subsystem \(p'\) as \(\Nsw_{p'}(s'-\varepsilon,s'+\varepsilon)\leq 0\). This limitation does not apply to the activation of subsystems \(p\in\P_S\) and/or switches from subsystems \(p\) to subsystems \(q\) with \((p,q)\in E_{-}(\P)\) because of the usage of ``at least (\(\geq\))'' relation with the total number of switches on any interval of time.}
    \end{remark}
     \begin{remark}
    \label{rem:sw-prop3}
    \rm{
        On every interval of time, \eqref{e:key_ineq5}-\eqref{e:key_ineq8} are inequality conditions and involve minimum (resp., maximum) number of activation of subsystem \(p\) (resp., switches between subsystems \(p\) and \(q\)), and \eqref{e:key_ineq4b}-\eqref{e:key_ineq4c} are equality conditions and involve the (exact) total number of activation of subsystem \(p\) (resp., switches between subsystems \(p\) and \(q\)) with a fixed \(\Nsw(s,t)\) obeying \eqref{e:key_ineq4a}. In other words, \eqref{e:key_ineq5}-\eqref{e:key_ineq8} are concerned with bounds and \eqref{e:key_ineq4b}-\eqref{e:key_ineq4c} are concerned with ``how'' close/far from the tight version of the bounds the frequencies could be chosen such that a given value of \(\Nsw(s,t)\) is obeyed.}
    \end{remark}
    
    \begin{theorem}
    \label{t:mainres}
    {\it
        Consider a family of systems \eqref{e:family}. Let the set of admissible switches between the subsystems, \(E(\P)\), and admissible minimum and maximum dwell times, \(\delta_p\) and \(\Delta_p\), on the subsystems \(p\in\P\), be given. Suppose that Assumptions \ref{assump:key1}-\ref{assump:key3} hold. Then the switched system \eqref{e:swsys} is input/output-to-state stable (IOSS) under every switching signal \(\sigma\in{\mathcal{S}}_{\mathcal{R}}'\).
        }
    \end{theorem}

    Our characterization of the class of stabilizing switching signals, \({\mathcal{S}}_{\mathcal{R}}'\), is the following: if (a) the subsystems admit linearly comparable Lyapunov-like functions (Assumptions \ref{assump:key1}-\ref{assump:key2}), and (b) the rates of decay (resp., growth) of Lyapunov-like functions corresponding to stable (resp., unstable) subsystems, the linear comparison factor between these functions, the admissible dwell times on the subsystems, and a set of scalars together satisfy a certain inequality (Assumption \ref{assump:key3}), then every admissible switching signal \(\sigma\in\mathcal{S}_{\mathcal{R}}\) for which the frequency of activation of various subsystems and frequency of switches between admissible pairs of subsystems, are governed by the set of scalars mentioned in (b) above (\'{a} la Definition \ref{d:adm-sw}) preserves IOSS of the switched system \eqref{e:swsys}. More specifically, a choice of \(\tilde{\mathcal{S}}_{\mathcal{R}}\) is \({\mathcal{S}}_{\mathcal{R}}'\). A proof of Theorem \ref{t:mainres} is presented in Section \ref{s:proof}. We will observe in the proof of Theorem \ref{t:mainres} that the functions \(\alpha\), \(\beta\), \(\chi_1\) and \(\chi_2\) can be chosen independent of \(\sigma\in\mathcal{S}_{\mathcal{R}}'\). Consequently, IOSS of the switched system \eqref{e:swsys} is \emph{uniform} over the elements of \(\mathcal{S}_{\mathcal{R}}'\) in the above sense. In the absence of outputs (resp. also inputs), the switched system \eqref{e:swsys} is ISS (resp., GAS) under the elements of \({\mathcal{S}}_{\mathcal{R}}'\).

    \begin{remark}
    \label{rem:compa}
    \rm{
        Notice that IOSS under arbitrary switching \cite{Mancilla2005} deals with the set of all switching signals, \(\Sw\), IOSS under average dwell time switching \cite{Muller2012} deals with a stabilizing subset \(\Sw'\) of \(\Sw\), and IOSS under given restrictions on the switching signals specifies \(\Sw_{\mathcal{R}}\subseteq\Sw\) a priori and deals with a stabilizing subset \(\Sw'_{\mathcal{R}}\) of \(\Sw_{\mathcal{R}}\). Our characterization of \(\Sw'_{\mathcal{R}}\) is more general to the existing characterization presented in \cite{def} in the following sense: (a) we allow (possibly different) admissible dwell times on the subsystems and (b) we do not restrict consecutive activation of unstable subsystems.
    }
    \end{remark}
    
     Given the Lyapunov-like functions, \(V_p\), \(p\in\P\) and their corresponding scalars, \(\lambda_p\), \(p\in\P\), \(\mu_{pq}\), \((p,q)\in E(\P)\), the scalars \(\rho^{S}_p\), \(p\in\P_S\), \(\rho^{U}_p\), \(p\in\P_U\), \(\rho^{+}_{pq}\), \((p,q)\in E_{+}(\P)\), \(\rho^{-}_{pq}\), \((p,q)\in E_{-}(\P)\) can be solved for as decision variables of a feasibility problem. Below we list a set of sufficient conditions under which Assumption \ref{assump:key3} is satisfied. 
    \begin{enumerate}[label = (\Roman*), leftmargin = 0.5cm]
        \item Suppose that \(V_p = V\) for all \(p\in\P\) and there exist \(\rho', \rho''\in[0,1[\) that satisfy \(\abs{\P_U}\rho'<1\), 
   \(\abs{\P_S}\rho''\leq1\), 
                \(\frac{\displaystyle{\rho'\Biggl(\sum_{p\in\P_U}\abs{\lambda_p}\Delta_p\Biggr)}}
                {\displaystyle{\rho''\Biggl(\sum_{p\in\P_S}\abs{\lambda_p}\delta_p\Biggr)}} < \frac{\delta_{\min}}{\Delta_{\max}}\).
            Then Assumption \ref{assump:key3} holds with \(\rho^{S}_p=\rho''\), \(p\in\P_S\), \(\rho^{U}_p=\rho'\), \(p\in\P_U\). Indeed, with \(\ln\mu_{pq}=\ln 1= 0\) for all \((p,q)\in E(\P)\), we have \(E_{+}(\P) = E_{-}(\P) = \emptyset\).
            
            \item Suppose that there exist  \(\rho', \rho''\in[0,1[\) that satisfy \(\abs{\P_U}\rho'<1\), \(\abs{E_{+}(\P)}\rho'\leq 1\),
   \(\abs{\P_S}\rho''\leq1\), \(\abs{E_{-}(\P)}\rho''\leq 1\),
                \(\frac{\rho'\Biggl(\displaystyle{\sum_{p\in\P_U}\abs{\lambda_u}\Delta_p+\sum_{(p,q)\in E_{+}(\P)}\abs{\ln\mu_{pq}}\Biggr)}}
                {\rho''\Biggl(\displaystyle{\sum_{p\in\P_S}\abs{\lambda_s}\delta_p+\sum_{(p,q)\in E_{-}(\P)}\abs{\ln\mu_{pq}}}\Biggr)} < \frac{\delta_{\min}}{\Delta_{\max}}\).
            Then Assumption \ref{assump:key3} holds with \(\rho^{S}_p = \rho''\), \(p\in\P_S\), \(\rho^{U}_p = \rho'\), \(p\in\P_U\), \(\rho^{+}_{pq} = \rho'\), \((p,q)\in E_{+}(\P)\), \(\rho^{-}_{pq} = \rho''\), \((p,q)\in E_{-}(\P)\).
            \item Suppose that \(\lambda_p=\lambda_s\) for all \(p\in\P_S\), \(\lambda_p=\lambda_u\) for all \(p\in\P_U\), \(\mu_{pq} = \mu^{+}\) for all \((p,q)\in E_{+}(\P)\), \(\mu_{pq} = \mu^{-}\) for all \((p,q)\in E_{-}(\P)\), \(\delta_p=\delta\), \(\Delta_p=\Delta\) for all \(p\in\P\), and there exist \(\rho'\), \(\rho''\in[0,1[\) that satisfy
               \(\abs{\P_S}\rho''\leq 1,\:\:\abs{\P_U}\rho' < 1\),
               \(\abs{E_{-}(\P)}\rho''\leq 1,\:\abs{E_{+}(\P)}\rho'\leq 1\), and
               \(\displaystyle{\frac{\rho'\Bigl(\abs{\lambda_u}\Delta\abs{\P_U}+\abs{\ln\mu^{+}}\abs{E_{+}(\P)}\Bigr)}
               {\rho''\Bigl(\abs{\lambda_s}\delta\abs{\P_S}+\abs{\ln\mu^{-}}\abs{E_{-}(\P)}\Bigr)} < \frac{\delta_{\min}}{\Delta_{\max}}}\).
            Then Assumption \ref{assump:key3} holds with \(\rho^{U}_{p}=\rho'\), \(p\in\P_U\), \(\rho^{+}_{pq} = \rho'\), \((p,q)\in E_{+}(\P)\), \(\rho^{S}_{p}=\rho''\), \(p\in\P_S\), \(\rho^{-}_{pq} = \rho''\), \((p,q)\in E_{-}(\P)\).
   \end{enumerate}
   We next provide a few examples where Assumption \ref{assump:key3} is satisfied.
   \begin{example}
   \label{ex:case1}
   \rm{
        Let \(\P_S = \{1\}\), \(\P_U = \{2\}\), \(E(\P) = \{(1,2),(2,1)\}\), \(\delta_1=3.4\), \(\Delta_1=4\), \(\delta_2 = 2\), \(\Delta_2=3.5\). Suppose that \(\lambda_1 = 1.25\), \(\lambda_2 = -0.5\), \(\mu_{12} = 1.25\), \(\mu_{21} = 0.8\). Then Assumption \ref{assump:key3} is satisfied with \(\rho^{S}_1=\rho^{U}_2=\rho^{+}_{12}=\rho^{-}_{21}=0.5\).
   }
   \end{example}
   
   \begin{example}
   \label{ex:case2}
   \rm{
        Let \(\P_S=\{1,2\}\), \(\P_U=\{3\}\), \(E(\P) = \{(1,2),(1,3)\), \((2,1),(3,1),(3,2)\}\), \(\delta_1=2\), \(\Delta_1=4\), \(\delta_2=2.5\), \(\Delta_2=3.5\), \(\delta_3=1.5\), \(\Delta_3=2\). Suppose that \(\lambda_1=1.5\), \(\lambda_2=1.25\), \(\lambda_3=-1.1\), \(\mu_{12}=\mu_{13}=1.25\), \(\mu_{21}=\mu_{31}=0.8\), \(\mu_{32}=0.6\). Then Assumption \ref{assump:key3} is satisfied with \(\rho^{S}_1=\rho^{U}_2=\rho^{U}_3 = \rho^{+}_{12}=\rho^{+}_{13}=\rho^{-}_{21}=\rho^{-}_{31}=\rho^{-}_{32}=0.3\).
   }
   \end{example}
   
   \begin{example}
   \label{ex:case3}
   \rm{
        Let \(\P_S=\{1,2,3,4,5\}\), \(\P_U=\{6,7,8,9,10\}\), \(E(\P) = \{(1,2),(1,3),(1,4),(1,5)\),\\\((2,1),(2,3),(2,4),(2,5)\), \((3,1),(3,2),(3,4),(3,5),(4,1),(4,2),(4,3),(4,5),
        (5,1)\), \((5,2)\),\\\((5,3),(5,4),(5,6),(6,1),(6,7),(7,2),(7,8),(8,3)\), \((8,9),(9,4),(9,10),(10,1),(10,2)\}\),\\\(\delta_1=\Delta_6 = 2\), \(\delta_2=\Delta_7 = 2.1\), \(\delta_3=\Delta_8 = 2.2\), \(\delta_4=\Delta_9 = 2.3\), \(\delta_5=\Delta_{10} = 2.4\), \(\delta_6 = 1\), \(\delta_7 = 1.1\), \(\delta_8 = 1.2\), \(\delta_9 = 1.3\), \(\delta_{10} = 1.4\), \(\Delta_1 = 3\), \(\Delta_2 = 3.1\), \(\Delta_3 = 3.2\), \(\Delta_4 = 3.3\), \(\Delta_5 = 3.4\). Suppose that \(\lambda_1=\lambda_2=\lambda_3=\lambda_4=\lambda_5=1.5\),
        \(\lambda_6=\lambda_7=\lambda_8=\lambda_9=\lambda_10=-0.75\), \(\mu_{12}=\mu_{13}=\mu_{14}=\mu_{15}=\mu_{21}=
        \mu_{23}=\mu_{24}=\mu_{25}=\mu_{31}=\mu_{32}=\mu_{34}=\mu_{35}=\mu_{41}=\mu_{42}=\mu_{43}=\mu_{45}=
        \mu_{51}=\mu_{52}=\mu_{53}=\mu_{54}=\mu_{67}=\mu_{78}=\mu_{89}=\mu_{9\:10} = 1\), \(\mu_{61}=\mu_{72}=\mu_{83}=\mu_{94}=\mu_{10\:1}=\mu_{10\:2}=0.8\), \(\mu_{56}=1.25\). Then Assumption \ref{assump:key3} is satisfied with \(\rho^{S}_1 = \rho^{S}_2 = \rho^{S}_3 = \rho^{S}_4 = \rho^{S}_5 =0.2\), \(\rho^{U}_6 = \rho^{U}_7 = \rho^{U}_8 = \rho^{U}_9 = \rho^{U}_{10} =0.1\),
        \(\rho^{-}_{61} = \rho^{-}_{72} = \rho^{-}_{83} = \rho^{-}_{94} = \rho^{-}_{10\:1} = \rho^{-}_{10\: 2} =0.1\), \(\rho^{+}_{56} = 0.1\).
   }
   \end{example}

   Co-designing multiple Lyapunov-like functions, \(V_p\), \(p\in\P\), their corresponding scalars, \(\lambda_p\), \(p\in\P\), \(\mu_{pq}\), \((p,q)\in E(\P)\) and  the scalars \(\rho^{S}_p\), \(p\in\P_S\), \(\rho^{U}_p\), \(p\in\P_U\), \(\rho^{+}_{pq}\), \((p,q)\in E_{+}(\P)\), \(\rho^{-}_{pq}\), \((p,q)\in E_{-}(\P)\) such that Assumption \ref{assump:key3} holds, is, however, a numerically difficult task. Further, since our stability condition is only sufficient, non-satisfaction of the Assumptions \ref{assump:key1}-\ref{assump:key3} does not imply the absence of stabilizing elements in \(\mathcal{S}_{\mathcal{R}}\).
\section{Numerical example}
\label{s:numex}
    We consider a family of systems \eqref{e:family} with 
    \begin{align*}
        f_p(x,v) = \pmat{a_{1p}x_1+b_{1p}sin(x_1-x_2)+c_{1p}v\\a_{2p}x_2+b_{2p}sin(x_2-x_1)+c_{2p}v},\:\:h_p(x) = x_1-x_2,
    \end{align*}
    where \(p\in\P = \{1,2,3,4\}\), \(a_{11}=a_{12}=a_{21}=a_{22}=-1\), \(a_{13}=a_{14}=a_{23}=a_{24} = 1\), 
    \(b_{11}=b_{12}=b_{13}=b_{14}=b_{23}=b_{24}=1\), \(b_{21}=b_{22}=0.8\), \(c_{11}=c_{13}=c_{21}=c_{23}=0\), 
    \(c_{12}=c_{14}=c_{22}=c_{24}=0.5\). We have \(\P_S=\{1,2\}\) and \(\P_U=\{3,4\}\). Suppose that the set of admissible switches is \(E(\P) = \{(1,2),(1,4),(2,1),(2,3),(3,4),(4,1)\}\) and the admissible dwell times are \(\delta_1=\delta_2=1\), \(\delta_3=1.5\), \(\delta_4=1.2\), \(\Delta_1=\Delta_3=2\), \(\Delta_2=\Delta_4=1.5\) units of time. We want to identify a class of switching signals that obeys the given constraints and preserves IOSS of the resulting switched system.
    \begin{remark}
    \label{rem:case_study}
    \rm{
        Notice that instead of considering IOSS under all switching signals \cite{Mancilla2005} or all switching signals satisfying a stabilizing average dwell time \cite{Muller2012}, we have restricted our attention to a pre-specified subset of the set of all switching signals. The existing literature on IOSS under restricted switching does not cater to the problem setting in the following sense: (a) application of the results in \cite{abc} gives \emph{a} stabilizing periodic switching signal instead of \emph{a class of} stabilizing switching signals, and (b) the results in \cite{def} does not apply as it is not allowed to switch from the unstable subsystem \(3\) to any of the stable subsystems.
        }
    \end{remark}
    We choose \(V_1(x)=V_2(x) = 0.5(x_1^2+1.25x_2^2)\), \(V_3(x)=V_4(x) = 0.5(x_1^2+x_2^2)\). Assumptions \ref{assump:key1} and \ref{assump:key2} are satisfied with \(\lambda_1=\lambda_2=1.75\), \(\lambda_3=\lambda_4=-2.1667\), \(\mu_{12}=\mu_{14}=\mu_{21}=\mu_{23}=\mu_{34}=1\), \(\mu_{41}=1.25\). We have
    that Assumption \ref{assump:key3} is satisfied with \(\rho^{S}_1=\rho^{S}_2=0.45\), \(\rho^{U}_{3}=\rho^{U}_{4}=\rho^{+}_{41}=0.1\). Indeed, 
    \begin{align*}
        &\rho^{S}_1+\rho^{S}_{2} = 0.9,\:\rho^{U}_3+\rho^{U}_4 = 0.2,\:\text{and}\\
        &-\frac{1}{\Delta_{\max}}\Biggl(\sum_{p\in\P_S}\abs{\lambda_p}\rho^{S}_{p}\delta_p+\sum_{(p,q)\in E_{-}(\P)}\abs{\ln\mu_{pq}}\rho^{-}_{pq}\Biggr)\nonumber\\
			&+\frac{1}{\delta_{\min}}\Biggl(\sum_{p\in\P_U}\abs{\lambda_p}\rho^{U}_{p}\Delta_p+\sum_{(p,q)\in E_{+}(\P)}\abs{\ln\mu_{pq}}\rho^{+}_{pq}\Biggr) = -0.0069 < 0.
    \end{align*}
    Consequently, the assertion of Theorem \ref{t:mainres} is applicable.
    
    \begin{figure}
        \centering
        \includegraphics[scale=0.4]{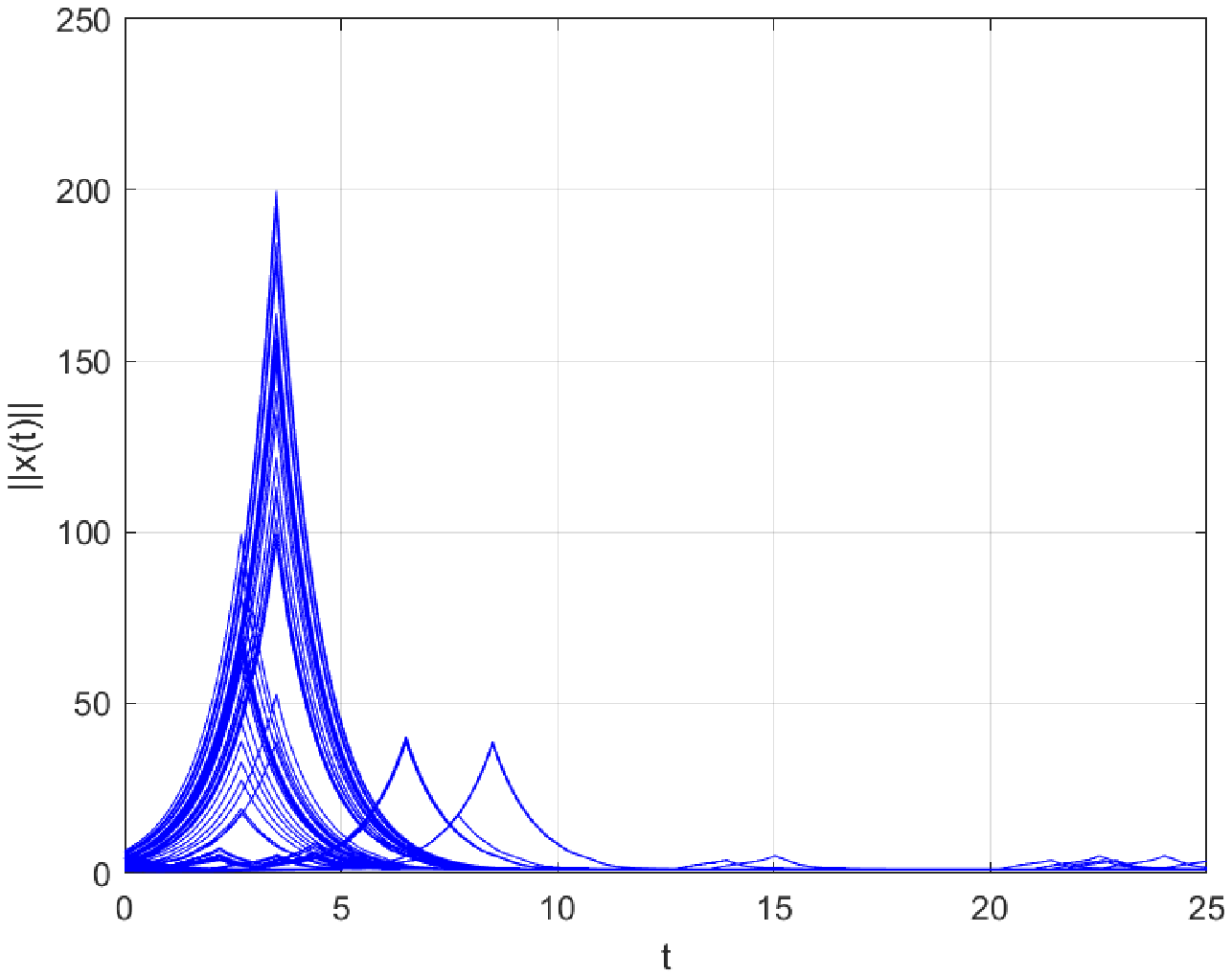}
        \caption{\(\bigl(\norm{x(t)}\bigr)_{t\geq 0}\) under \(\sigma\in\Sw'_{\mathcal{R}}\)}\label{fig:ioss}
    \end{figure}
    
    We fix  \(\tilde{\rho}^{U}_{3}=\tilde{\rho}^{U}_{4}=\tilde{\rho}^{+}_{41} = 1\). Notice that the number of switches for \((1,2),(1,4),(2,1),(2,3)\),\\\((3,4)\) can be chosen freely to match the total number of switches. We generate \(10\) elements of the set \(\Sw'_{\mathcal{R}}\) and corresponding to each element, we choose 10 different \(x_0\in[-5,5]^2\) uniformly at random. The input is chosen uniformly at random from the range \(v\in[-0.5,0.5]\). The resulting processes \(\bigl(\norm{x(t)}\bigr)_{t\geq 0}\) are illustrated till \(t=25\) units of time in Fig. \ref{fig:ioss}. IOSS of the switched system \eqref{e:swsys} is observed. 
\section{Concluding remarks}
\label{s:concln}
	In this paper we identified a class of switching signals that guarantees IOSS of a continuous-time switched nonlinear system under restrictions on admissible switches between the subsystems and admissible dwell times on the subsystems. We identify the design of state-norm estimators (\`a la \cite[\S 4]{Muller2012}) under the proposed class of stabilizing switching signals as a future research direction.
\section{Proof of Theorem \ref{t:mainres}}
\label{s:proof}
		Fix \(\sigma\in\mathcal{S}_{\mathcal{R}}'\) and \(t > 0\). Let \(0=:\tau_0<\tau_1<\cdots<\tau_{\Nsw(0,t)}\) be the switching instants before (and including) \(t\). In view of \eqref{e:key_ineq2}, we have that
		\begin{align*}
			V_{\sigma(t)}(x(t))&\leq \exp\Bigl(-\lambda_{\sigma(\tau_{\Nsw(0,t)})}(t-\tau_{\Nsw(0,t)})\Bigr)V_{\sigma(t)}\Bigl(x(\tau_{\Nsw(0,t)})\Bigr)\\
			&\hspace*{-1.5cm}+\Bigl(\gamma_1\bigl(\norm{v}_{[0,t]}\bigr)+\gamma_2\bigl(\norm{y}_{[0,t]}\bigr)\Bigr)
			\times
			\int_{\tau_{\Nsw(0,t)}}^{t}\exp\Bigl(-\lambda_{\sigma(\tau_{\Nsw(0,t)})}(t-s)\Bigr)ds.
		\end{align*}
		Applying \eqref{e:key_ineq3} and iterating the above, we obtain the estimate
		\begin{align}
		\label{e:pf_step1}
			V_{\sigma(t)}(x(t))&\leq\psi_1(t)V_{\sigma(0)}(x(0))
			 + \Bigl(\gamma_1\bigl(\norm{v}_{[0,t]}\bigr)+\gamma_2\bigl(\norm{y}_{[0,t]}\bigr)\Bigr)\psi_2(t),
		\end{align}
		where
		\begin{align}
		\label{e:psi1_defn}
			\psi_1(t) &:= \exp\Biggl(-\sum_{\substack{{i=0}\\{\tau_{\Nsw(0,t)+1}:=t}}}^{\Nsw(0,t)}\lambda_{\sigma(\tau_i)}(\tau_{i+1}-\tau_{i})
			+\sum_{i=0}^{\Nsw(0,t)-1}\ln\mu_{\sigma(\tau_i)\sigma(\tau_{i+1})}\Biggr),
		\end{align}
		and
		\begin{align}
		\label{e:psi2_defn}
			\psi_2(t) &:= \sum_{\substack{{i=0}\\{\tau_{\Nsw(0,t)+1}:=t}}}^{\Nsw(0,t)}\Biggl(\exp\Biggl(-\sum_{\substack{{j=i+1}\\{\tau_{\Nsw(0,t)+1}:=t}}}^{\Nsw(0,t)}\lambda_{\sigma(\tau_j)}(\tau_{j+1}-\tau_{j})\nn\\&\hspace*{-0.7cm}+\sum_{j=i+1}^{\Nsw(0,t)-1}\ln\mu_{\sigma(\tau_j)\sigma(\tau_{j+1})}\Biggr)
			\times
			\frac{1}{\lambda_{\sigma(\tau_i)}}\Bigl(1-\exp\bigl(-\lambda_{\sigma(\tau_i)}(\tau_{i+1}-\tau_{i})\bigr)\Bigr)\Biggr).
		\end{align}
		
		In view of \eqref{e:key_ineq1}, we rewrite the estimate \eqref{e:pf_step1} as
		\begin{align*}
			\underline{\alpha}(\norm{x(t)})&\leq\psi_1(t)\overline{\alpha}(\norm{x(0)})
			+\bigl(\gamma_1(\norm{v}_{[0,t]})+\gamma_2(\norm{y}_{[0,t]})\bigr)\psi_2(t).
		\end{align*}
		By Definition \ref{d:ioss-cont}, for IOSS of \eqref{e:swsys}, we need to show that
			 i) \(\overline{\alpha}(\star)\psi_1(\cdot)\) can be bounded above by a class \(\KL\) function, and
			ii) \(\psi_2(\cdot)\) is bounded above by a constant.
		
		We first verify i). Towards this end, we already see that \(\overline{\alpha}\in\Kinfty\) from Assumption \ref{assump:key1}. Therefore, it remains to show that \(\psi_1(\cdot)\) is bounded above by a function in class \(\Lf\).
		
		Now, the argument of the exponential in \(\psi_1(t)\) is equal to
			\(\displaystyle{-\sum_{p\in\P_S}\abs{\lambda_p}\T_p(0,t)}\)\\\(\displaystyle{-\sum_{(p,q)\in E_{-}(\P)}\abs{\ln\mu_{pq}}\Nsw_{pq}(0,t)}\)
			\(\displaystyle{+\sum_{p\in\P_U}\abs{\lambda_p}\T_p(0,t)+\sum_{(p,q)\in E_{+}(\P)}\abs{\ln\mu_{pq}}\Nsw_{pq}(0,t)}\).
		In view of 
		 properties of \(\sigma\in\mathcal{S}_{\mathcal{R}}'\) and the fact that for \(a\in\R\), \(a-1\leq\lfloor a\rfloor\leq a\), the above expression is bounded above by
		\begin{align*}
			&\displaystyle{{c'} + \Biggl(-\sum_{p\in\P_S}\abs{\lambda_p}\rho^{\mathrm{S}}_{p}\frac{\delta_p}{\Delta_{\max}}}\displaystyle{-\sum_{(p,q)\in E_{-}(\P)}\abs{\ln\mu_{pq}}\rho^{-}_{pq}\frac{1}{\Delta_{\max}}}\\
			&\quad\quad+\sum_{p\in\P_U}\abs{\lambda_p}\rho^{\mathrm{U}}_{p}\frac{\Delta_p}{\delta_{\min}}
			\displaystyle{+\sum_{(p,q)\in E_{+}(\P)}\abs{\ln\mu_{pq}}\rho^{\mathrm{+}}_{pq}\frac{1}{\delta_{\min}}\Biggr)t},
		\end{align*}
		where \({c'}\) is \(\displaystyle{\sum_{p\in\P_U}\abs{\lambda_p}\tilde{\rho}^{\mathrm{U}}_{p}\Delta_p+
        \sum_{(p,q)\in E_{+}(\P)}\abs{\ln\mu_{p}}\tilde{\rho}^{+}_{pq}}\)
        \(\displaystyle{+\sum_{p\in\P_S}\abs{\lambda_p}\delta_p}\)\(\displaystyle{+\sum_{(p,q)\in E_{-}(\P)}\abs{\ln\mu_{pq}}}\)\\
        \(\displaystyle{+\Biggl(\sum_{p\in\P_S}\abs{\lambda_p}\rho^{\mathrm{S}}_{p}\delta_p
        +\sum_{(p,q)\in E_{-}(\P)}\abs{\ln\mu_{pq}}\rho^{-}_{pq}\Biggr)}\).
		
		By \eqref{e:key_ineq4}, we have that the above expression is at most equal to \(c_1-c_2 t\) for some \(c_1\geq {c'}\) and \(c_2>0\). Therefore, from \eqref{e:psi1_defn}, we obtain that \(\psi_1(t)\leq\exp\bigl(c_1-c_2 t\bigr)\).
		Notice that the right-hand side of the above inequality decreases as \(t\) increases and tends to \(0\) as \(t\to+\infty\). Consequently, i) holds.
		
		We now verify ii). We have that the function \(\psi_{2}(t)\) is bounded above by
        \begin{align}
        \label{e:pf1_step2}
            &\sum_{p\in\P_{S}}\frac{1}{\abs{\lambda_{p}}}\sum_{\substack{i=0\\\sigma(\tau_i)=p\\\tau_{\Nsw(0,t)+1}:=t}}^{\Nsw(0,t)}\Biggl(\exp\Biggl(-\sum_{k\in\P_{S}}\abs{\lambda_{k}}\T_{k}(\tau_{i+1},t)
            +\sum_{k\in\P_{U}}\abs{\lambda_{k}}\T_{k}(\tau_{i+1},t)\\
            &\quad\quad-\sum_{(k,\ell)\in E_{-}(\P)}\abs{\ln\mu_{k\ell}}\Nsw_{k\ell}(\tau_{i+1},t)
            +\sum_{(k,\ell)\in E_{+}(\P)}(\ln\mu_{k\ell})\Nsw_{k\ell}(\tau_{i+1},t)\Biggr)\Biggr)\nn\\
            -&\sum_{p\in\P_{U}}\frac{1}{\abs{\lambda_{p}}}\sum_{\substack{i=0\\\sigma(\tau_{i})=p\\\tau_{\Nsw(0,t)+1}:=t}}^{\N(0,t)}\Biggl(\exp\Biggl(-\sum_{k\in\P_{S}}\abs{\lambda_{k}}\T_{k}(\tau_{i},t)
            +\sum_{k\in\P_{U}}\abs{\lambda_{k}}\T_{k}(\tau_{i},t)\\
            &\quad\quad-\sum_{(k,\ell)\in E_{-}(\P)}\abs{\ln\mu_{k\ell}}\Nsw_{k\ell}(\tau_{i},t)
            +\sum_{(k,\ell)\in E_{+}(\P)}(\ln\mu_{k\ell})\Nsw_{k\ell}(\tau_{i},t)\Biggr)\Biggr).
          \end{align}
         Now, by an application of a similar set of arguments as employed in the verification of i), we obtain that
         \begin{align*}
         	&\hspace*{-0.2cm}\displaystyle{-\sum_{k\in\P_{S}}\abs{\lambda_{k}}\T_{k}(\tau_{i},t)
            +\sum_{k\in\P_{U}}\abs{\lambda_{k}}\T_{k}(\tau_{i},t)}\\
            &\quad\quad\displaystyle{-\sum_{(k,\ell)\in E_{-}(\P)}\abs{\ln\mu_{k\ell}}\Nsw_{k\ell}(\tau_{i},t)
            +\sum_{(k,\ell)\in E_{+}(\P)}(\ln\mu_{k\ell})\Nsw_{k\ell}(\tau_{i},t)}\\
            \leq&\displaystyle{c''+\Biggl(-\sum_{k\in\P_{S}}\abs{\lambda_{k}}\rho^{\mathrm{S}}_{k}\frac{\delta_k}{\Delta_{\max}}}-\sum_{(k,\ell)\in E_{-}(\P)}\abs{\ln\mu_{k\ell}}\rho^{-}_{k\ell}\frac{1}{\Delta_{\max}}\\
              &\quad\quad\displaystyle{+\sum_{k\in\P_{U}}\abs{\lambda_{k}}\rho^{\mathrm{U}}_{k}\frac{\Delta_k}{\delta_{\min}}}
            \quad\quad\displaystyle{+\sum_{(k,\ell)\in E_{+}(\P)}(\ln\mu_{k\ell})\rho^{+}_{k\ell}\frac{1}{\delta_{\min}}
         \Biggr)(t-\tau_i)},
           \end{align*}
           where \(c''\) is \(\displaystyle{\sum_{k\in\P_U}\abs{\lambda_k}\tilde{\rho}^{\mathrm{U}}_{k}\Delta_k+
           \sum_{(k,\ell)\in E_{+}(\P)}\abs{\ln\mu_{k\ell}}\tilde{\rho}^{+}_{k\ell}}\)
           \(\displaystyle{+\sum_{k\in\P_S}\abs{\lambda_k}\delta_k}\)
           \(\displaystyle{+\sum_{(k,\ell)\in E_{-}(\P)}\abs{\ln\mu_{k\ell}}}\)\\
           \(\displaystyle{+\sum_{k\in\P_S}\abs{\lambda_k}\rho^{\mathrm{S}}_{k}\delta_k+\sum_{(k,\ell)\in E_{-}(\P)}\abs{\ln\mu_{k\ell}}\rho^{-}_{k\ell}}\).
          
           In view of \eqref{e:key_ineq4}, the above expression is at most equal to \(\overline{c}_1-\overline{c}_2(t-\tau_i)\) for some \(\overline{c}_1\geq c''\) and \(\overline{c}_2>0\). Similarly, the expression
         	\(\displaystyle{-\sum_{k\in\P_{S}}\abs{\lambda_{k}}\T_{k}(\tau_{i+1},t)-\sum_{(k,\ell)\in E_{-}(\P)}\abs{\ln\mu_{k\ell}}\Nsw_{k\ell}(\tau_{i+1},t)}\)
            \\
            \(\displaystyle{+\sum_{k\in\P_{U}}\abs{\lambda_{k}}\T_{k}(\tau_{i+1},t)}
            \displaystyle{+\sum_{(k,\ell)\in E_{+}(\P)}(\ln\mu_{k\ell})\Nsw_{k\ell}(\tau_{i+1},t)}\)
            is at most equal to \(\tilde{c}_1-\tilde{c}_2(t-\tau_{i+1})\) for some \(\tilde{c}_1\geq c''\) and \(\tilde{c}_2 > 0\).

            Now, from \eqref{e:pf1_step2}, we have that 
              \begin{align*}
              	\psi_2(t) &\leq \sum_{p\in\P_S}\frac{1}{\abs{\lambda_p}}\sum_{\substack{i=0\\\sigma(\tau_i)=p\\\tau_{\Nsw(0,t)+1:=t}}}^{\Nsw(0,t)}\exp\bigl(\tilde{c}_1-\tilde{c}_2(t-\tau_{i+1})\bigr)\\ &\hspace*{2cm}+ \sum_{p\in\P_U}\frac{1}{\abs{\lambda_p}}\sum_{\substack{i=0\\\sigma(\tau_i)=p\\\tau_{\Nsw(0,t)+1:=t}}}^{\Nsw(0,t)}\exp\bigl(\overline{c}_1-\overline{c}_2(t-\tau_{i})\bigr)\\
	&\leq \sum_{p\in\P_S}\frac{1}{\abs{\lambda_p}}\sum_{i=0}^{\Nsw(0,t)}\exp\bigl(\tilde{c}_1-\tilde{c}_2(t-\tau_{i+1})\bigr)\\
	&\hspace*{2cm}+ \sum_{p\in\P_U}\frac{1}{\abs{\lambda_p}}\sum_{i=0}^{\Nsw(0,t)}\exp\bigl(\overline{c}_1-\overline{c}_2(t-\tau_{i})\bigr).
        \end{align*}
       We have
        \begin{align}
         \label{e:pf_step7}
            \sum_{i=0}^{\Nsw(0,t)}\exp\biggl(-\overline{c}_2(t-\tau_{i})\biggr) &\leq \sum_{i=0}^{\big\lfloor\frac{t}{\delta_{\min}}\big\rfloor}\exp\biggl(-\overline{c}_2(t-\tau_{i})\biggr)\nonumber\\
            =&\:\exp\biggl(-\overline{c}_2(t-0)\biggr)+\cdots+
            \exp\biggl(-\overline{c}_2(t-\tau_{\lfloor\frac{t}{\delta_{\min}}\rfloor})\biggr)\biggr)\nonumber\\
            \leq&\:\exp\bigl(-\overline{c}_2t\bigr)+\exp\biggl(-\overline{c}_2\bigg\lfloor\frac{t}{\Delta_{\max}}\bigg\rfloor\delta_{\min}\biggr)+
          \cdots+
          +\exp(-\overline{c}_2\delta_{\min})\nonumber\\
            \leq&\:1 + \exp\bigl(-\overline{c}_2\delta_{\min}\bigr)\frac{1-\bigl(\exp\bigl(-\overline{c}_2\delta_{\min}\bigr)\bigr)^{\big\lfloor\frac{t}{\Delta_{\max}}
            \big\rfloor+1}}{1-\exp\bigl(-\overline{c}_2\delta_{\min}\bigr)}\nonumber\\
            \leq&\:1+\frac{1}{\exp\bigl(\overline{c}_2\delta_{\min}\bigr)-1}.
         \end{align}
         Applying the set of arguments as above, we also obtain that
         \begin{align}
         \label{e:pf_step8}
            \sum_{i=0}^{\Nsw(0,t)}\exp\biggl(-\tilde{c}_2\bigl(t-\tau_{i+1}\bigr)\biggr)\leq \frac{1}{\exp\bigl(\tilde{c}_2\delta_{\min}\bigr)-1}.
         \end{align}
         In view of \eqref{e:pf_step7}-\eqref{e:pf_step8}, we arrive at ii). Indeed,
         \begin{align*}
            \psi_2(t)\leq \overline{\psi_2} &= \sum_{p\in\P_S}\frac{1}{\abs{\lambda_p}}\exp(\tilde{c}_1)\frac{1}{\exp(\tilde{c}_2\delta)-1}
            +\sum_{p\in\P_U}\frac{1}{\abs{\lambda_p}}\exp(\overline{c}_1)\biggl(1+\frac{1}{\exp(\overline{c}_2\delta)-1}\biggr).
         \end{align*}

	We conclude IOSS of the switched system \eqref{e:swsys} under \(\sigma\). Recall that \(\sigma\in\mathcal{S}_{\mathcal{R}}'\) was chosen arbitrarily. Consequently, the assertion of Theorem \ref{t:mainres} follows. In particular, condition \eqref{e:ioss-cont} holds with \(\alpha(r) = r\), \(\beta(r,s) = \overline{\alpha}(r)\exp\bigl(c_1-c_2s\bigr)\), \(\chi_1(r) = \gamma_1(r)\overline{\psi}_2\) and \(\chi_2(r) = \gamma_2(r)\overline{\psi}_{2}\).
        \hspace*{3.5cm}\qed

\end{document}